\ifx\shlhetal\undefinedcontrolsequenc\let\shlhetal\relax\fi
%%%%%%%%%%%%%%%%%%%%%%%%%%%%%%%%%%%%%%%%%%%%%%%%%%%%%%%%%%%%%%%%
%%%%%%% final changes and improvements were done by Saharon on
%%%%%%% February 12-14, 1995
%%%%%%% corrections - galley proof in March 1995
%%%%%%%%%%%%%%%%%%%%%%%%%%%%%%%%%%%%%%%%%%%%%%%%%%%%%%%%%%%%%%%%
% more small corrections. goldstrn, Aug 22-25
% Fact 4 added. goldstrn, Aug 14, 94.
% modification to references by  goldstrn, Aug 14, 94.

% From leonhard@math.rutgers.edu Fri Jul  8 17:22:56 1994
% Date: Fri, 8 Jul 94 10:19:03 EDT
% From: Alice Leonhardt <leonhard@math.rutgers.edu>
% To: Saharon Shelah's office <shlhetal@math.huji.ac.il>

\def\extend{{{}^\frown}\!}
\def\IT{\fam\itfam}
\input amstex

\def\restriction{\mathord{\upharpoonright}}

%%%\magnification=\magstep 1
\documentstyle {amsppt}

%%% \input bib4plain
%
% file bib4plain.tex starts here 
%
%  This file should be inputted if you want to use 
%  bibtex fom within plain TeX. 
      % Not really need for standard
       % bibtex files, but these commands
\def\renewcommand{\newcommand}	       % are used in our literal-unsrt.bst
\edef\cite{\the\catcode`@}%
\catcode`@ = 11
\let\@oldatcatcode = \cite
\chardef\@letter = 11
\chardef\@other = 12
%
%
% Next come some things that will be useful later.
%
% Make an outer definition into an inner one (due to Chris Thompson).
% The arguments should be the control sequence to be defined, and the
% new of the \outer control sequence, as characters; the control
% sequence #1 is defined to be just the same as \csname#2\endcsname, but
% not \outer.  For example, \@innerdef\innernewcount{newcount} would
% define \innernewcount to be a non-outer version of \newcount.
%
\def\@innerdef#1#2{\edef#1{\expandafter\noexpand\csname #2\endcsname}}%
%
% We use \@innerdef to make some of our allocations local, because
% Eplain includes our code inside a conditional.  We put @'s in the
% names to minimize the (already small) chance of conflicts.
%
\@innerdef\@innernewcount{newcount}%
\@innerdef\@innernewdimen{newdimen}%
\@innerdef\@innernewif{newif}%
\@innerdef\@innernewwrite{newwrite}%
%
%
% Swallow one parameter.
%
\def\@gobble#1{}%
%
%
% Use TeX 3.0's \inputlineno to get the line number, for better error
% messages, but if we're using an old version of TeX, don't do anything.
%
\ifx\inputlineno\@undefined
   \let\@linenumber = \empty % Pre-3.0.
\else
   \def\@linenumber{\the\inputlineno:\space}%
\fi
%
%
% The following macro \@futurenonspacelet (from the TeXbook) behaves
% essentially like \futurelet except that it discards any implicit or
% explicit space tokens that intervene before a nonspace is scanned:
%
\def\@futurenonspacelet#1{\def\cs{#1}%
   \afterassignment\@stepone\let\@nexttoken=
}%
\begingroup % The grouping here avoids stepping on an outside use of `\\'.
\def\\{\global\let\@stoken= }%
\\ % now \@stoken is a space token (\\ is a control symbol, so that
   % space after it is seen).
\endgroup
\def\@stepone{\expandafter\futurelet\cs\@steptwo}%
\def\@steptwo{\expandafter\ifx\cs\@stoken\let\@@next=\@stepthree
   \else\let\@@next=\@nexttoken\fi \@@next}%
\def\@stepthree{\afterassignment\@stepone\let\@@next= }%
%
%
% \@getoptionalarg\CS gets an optional argument from the input, enclosed
% in brackets, then expands \CS.  We set \@optionalarg to \empty if we
% don't find one, otherwise to the text of the argument.  This assumes
% the brackets don't have a funny category code.
%
\def\@getoptionalarg#1{%
   \let\@optionaltemp = #1%
   \let\@optionalnext = \relax
   \@futurenonspacelet\@optionalnext\@bracketcheck
}%
%
% The \expandafter's in this macro let us avoid the use of \aftergroup,
% which is somewhat more expensive.
%
\def\@bracketcheck{%
   \ifx [\@optionalnext
      \expandafter\@@getoptionalarg
   \else
      \let\@optionalarg = \empty
      % We can't do the \temp after the \fi, because then the \temp gets
      % in the way of reading the optional argument from the input, if
      % we do have one.
      \expandafter\@optionaltemp
   \fi
}%
\def\@@getoptionalarg[#1]{%
   \def\@optionalarg{#1}%
   \@optionaltemp
}%
%
%
% From LaTeX.
%
\def\@nnil{\@nil}%
\def\@fornoop#1\@@#2#3{}%
\def\@for#1:=#2\do#3{%
   \edef\@fortmp{#2}%
   \ifx\@fortmp\empty \else
      \expandafter\@forloop#2,\@nil,\@nil\@@#1{#3}%
   \fi
}%
\def\@forloop#1,#2,#3\@@#4#5{\def#4{#1}\ifx #4\@nnil \else
       #5\def#4{#2}\ifx #4\@nnil \else#5\@iforloop #3\@@#4{#5}\fi\fi
}%
\def\@iforloop#1,#2\@@#3#4{\def#3{#1}\ifx #3\@nnil
       \let\@nextwhile=\@fornoop \else
      #4\relax\let\@nextwhile=\@iforloop\fi\@nextwhile#2\@@#3{#4}%
}%
%
%
% This macro tests if a file \jobname.#1 exists, and sets \if@fileexists
% appropriately.  If an optional argument is given, it is used as the
% root part of the filename instead of \jobname.
%
\@innernewif\if@fileexists
\def\@testfileexistence{\@getoptionalarg\@finishtestfileexistence}%
\def\@finishtestfileexistence#1{%
   \begingroup
      \def\extension{#1}%
      \immediate\openin0 =
         \ifx\@optionalarg\empty\jobname\else\@optionalarg\fi
         \ifx\extension\empty \else .#1\fi
         \space
      \ifeof 0
         \global\@fileexistsfalse
      \else
         \global\@fileexiststrue
      \fi
      \immediate\closein0
   \endgroup
}%
%
%
%% [[[start of BibTeX-specific stuff]]]
%
% Now come the four main LaTeX commands and their associated .aux
% commands.  Just as in LaTeX, \bibliographystyle defines the BibTeX
% style name (.bst file, that is), and \bibliography defines the
% database (.bib) file(s).  The corresponding .aux-file commands are
% \bibstyle and \bibdata, which are there only for BibTeX's (but not
% LaTeX's) use.
%
\def\bibliographystyle#1{%
   \@readauxfile
   \@writeaux{\string\bibstyle{#1}}%
}%
\let\bibstyle = \@gobble
%
% As well as writing the \bibdata command to tell BibTeX which .bib
% files to read, we read the .bbl file that BibTeX (or a person,
% conceivably) has produced.  We use \bblfilebasename as the root of the
% filename to read; this defaults to \jobname.
%
\let\bblfilebasename = \jobname
\def\bibliography#1{%
   \@readauxfile
   \@writeaux{\string\bibdata{#1}}%
   \@testfileexistence[\bblfilebasename]{bbl}%
   \if@fileexists
      % We just output a non-discardable item (the `whatsit' with the
      % \bibdata command).  This means that the glue that will be
      % inserted next (\parskip or \baselineskip, most likely) will be a
      % legal breakpoint.  Most likely, this is after some kind of
      % heading, where we don't want to allow a page break.  So:
      \nobreak
      \@readbblfile
   \fi
}%
\let\bibdata = \@gobble
%
% The \nocite{label,label,...} command writes its argument to \@auxfile,
% unless instructed not to, but produces no text in the document.  Both
% the \nocite and \cite commands produce \citation commands in the .aux file.
%
\def\nocite#1{%
   \@readauxfile
   \@writeaux{\string\citation{#1}}%
}%
\@innernewif\if@notfirstcitation
%
% \cite[note]{label,label,...} produces the citations for the labels as
% well.  If the optional argument `note' is present, it's added after
% the labels.  Since \cite calls \nocite to do its .aux-file writing,
% \cite doesn't need to call \@readauxfile (\nocite does).
%
\def\cite{\@getoptionalarg\@cite}%
%
% Typeset the citations for the labels in #1, followed by the note, if
% it exists.  To change the citation's format in the text, redefine one
% or more `\print...' macros, whose defaults appear later in this file.
%
\def\@cite#1{%
   % Remember the optional argument, in case one of the macros we call
   % below ends up looking for an optional argument itself.  For
   % example, if a \cite[note] triggers reading the .aux file, then the
   % [note] would be clobbered, since \@testfileexistence looks for an
   % optional arg.
   \let\@citenotetext = \@optionalarg
   % Start printing the text, beginning with a left bracket by default.
   \printcitestart
   % It's complicated, but because \nocite puts a `whatsit' onto the list,
   % \nocite should follow \printcitestart.  It's conceivable, but very
   % unlikely, that this `whatsit' will cause a problem (glue that doesn't
   % disappear when you want it to is the most likely symptom), requiring
   % a change either to \printcitestart or to the label that the .bst file
   % produces.
   \nocite{#1}%
   \@notfirstcitationfalse
   \@for \@citation :=#1\do
   {%
      \expandafter\@onecitation\@citation\@@
   }%
   \ifx\empty\@citenotetext\else
      \printcitenote{\@citenotetext}%
   \fi
   \printcitefinish
}%
\def\@onecitation#1\@@{%
   \if@notfirstcitation
      \printbetweencitations
   \fi
   \expandafter \ifx \csname\@citelabel{#1}\endcsname \relax
      \if@citewarning
         \message{\@linenumber Undefined citation `#1'.}%
      \fi
      % Give it a dummy definition:
      \expandafter\gdef\csname\@citelabel{#1}\endcsname{%
         {\tt
            \escapechar = -1
            \nobreak\hskip0pt
            \expandafter\string\csname#1\endcsname
            \nobreak\hskip0pt
         }%
      }%
   \fi
   % Now produce the text, whether it was undefined or not.
   \csname\@citelabel{#1}\endcsname
   \@notfirstcitationtrue
}%
%
% Given a label `foo', the macro `\b@foo' is supposed to
% hold the text that should be produced.
%
\def\@citelabel#1{b@#1}%
%
% So, how does a citation label get defined?  When we read the .bbl file
% (below), a \bibitem writes out a \@citedef command.  And when we read
% the \@citedef, we define \@citelabel{#1}, where #1 is the user's
% label.
%
\def\@citedef#1#2{\expandafter\gdef\csname\@citelabel{#1}\endcsname{#2}}%
%
%
% Reading the .bbl file also produces the typeset bibliography.  Please
% notice, however, that we do not produce the title for the references
% (e.g., `References'), as LaTeX does.  The formatting and spacing of
% that title, whether it should go into the headline, and so on, are all
% things determined by your format.  We cannot know those things in
% advance.  If you wish, you can define \bblhook to produce the title.
% Or just do it before the \bibliography command.
%
\def\@readbblfile{%
   % Define a counter to tell us which item number we are on, unless
   % we've already defined it (because the document has more than one
   % bibliography).
   \ifx\@itemnum\@undefined
      \@innernewcount\@itemnum
   \fi
   \begingroup
      \def\begin##1##2{%
         % ##1 is just `thebibliography'.
         % ##2 is the widest label.
         % We set (new dimen) \biblabelwidth based on the widest label
         \setbox0 = \hbox{\biblabelcontents{##2}}%
         \biblabelwidth = \wd0
      }%
      \def\end##1{}% ##1 is `thebibliography' again.
      %
      % Here we have two possibilities:
      % \bibitem[typesetlabel]{citationlabel}
      % \bibitem{citationlabel}
      % If we have the second of these, the citations are numbered, starting
      % from one; we use our own count register \@itemnum for this.
      %
      \@itemnum = 0
      \def\bibitem{\@getoptionalarg\@bibitem}%
      \def\@bibitem{%
         \ifx\@optionalarg\empty
            \expandafter\@numberedbibitem
         \else
            \expandafter\@alphabibitem
         \fi
      }%
      \def\@alphabibitem##1{%
         % Need \xdef here for various reasons.
         \expandafter \xdef\csname\@citelabel{##1}\endcsname {\@optionalarg}%
         % Left-justify alpha labels, unless \biblabel{pre,post}contents
         % are already defined.
         \ifx\biblabelprecontents\@undefined
            \let\biblabelprecontents = \relax
         \fi
         \ifx\biblabelpostcontents\@undefined
            \let\biblabelpostcontents = \hss
         \fi
         \@finishbibitem{##1}%
      }%
      \def\@numberedbibitem##1{%
         \advance\@itemnum by 1
         \expandafter \xdef\csname\@citelabel{##1}\endcsname{\number\@itemnum}%
         % Right-justify numeric labels, unless \biblabel{pre,post}contents
         % are already defined.
         \ifx\biblabelprecontents\@undefined
            \let\biblabelprecontents = \hss
         \fi
         \ifx\biblabelpostcontents\@undefined
            \let\biblabelpostcontents = \relax
         \fi
         \@finishbibitem{##1}%
      }%
      \def\@finishbibitem##1{%
         \biblabelprint{\csname\@citelabel{##1}\endcsname}%
         \@writeaux{\string\@citedef{##1}{\csname\@citelabel{##1}\endcsname}}%
         \ignorespaces
      }%
      %
      % Do the printing (we're producing the bibliography, remember).
      %
      \let\em = \bblem
      \let\newblock = \bblnewblock
      \let\sc = \bblsc
      % Punctuation won't affect spacing;
      \frenchspacing
      % the penalties below are from LaTeX's [article,book,report].sty;
      \clubpenalty = 4000 \widowpenalty = 4000
      % the next two values come from LaTeX's \sloppy command;
      \tolerance = 10000 \hfuzz = .5pt
      \everypar = {\hangindent = \biblabelwidth
                      \advance\hangindent by \biblabelextraspace}%
      \bblrm
      % the \parskip is a guess at what looks good;
      \parskip = 1.5ex plus .5ex minus .5ex
      % and the space between label and text comes from LaTeX's \labelsep.
      \biblabelextraspace = .5em
      \bblhook
      \input \bblfilebasename.bbl
    \endgroup
}%
%
% The widest label's width is useful for redefining \biblabelprint;
% you redefine \biblabelwidth, in effect, by redefining the
% \biblabelcontents macro that appears below.  And \biblabelextraspace,
% which is redefinable inside \bblhook, is added to \biblabelwidth to
% determine the amount of hanging indentation.
%
\@innernewdimen\biblabelwidth
\@innernewdimen\biblabelextraspace
%
% Now come the main macros that are related to the printing of the
% bibliography.  Since you might want to redefine them, they are given
% default definitions outside of \@readbblfile.
%
% The first one controls the printing of a bibliography entry's label.
% If you change it, make sure that it starts with something like
% \noindent or \indent or \leavevmode that puts TeX into horizontal mode
% (even if the label itself is empty); otherwise, the hanging
% indentation will get messed up in certain circumstances.
%
\def\biblabelprint#1{%
   \noindent
   \hbox to \biblabelwidth{%
      \biblabelprecontents
      \biblabelcontents{#1}%
      \biblabelpostcontents
   }%
   \kern\biblabelextraspace
}%
%
% If you are using numeric labels, and you want them left-justified
% (numeric labels by default are right-justified), do something like:
%     \def\biblabelprecontents{\relax}
%     \def\biblabelpostcontents{\hss}
%
% By default the labels are typeset in \bblrm, and enclosed in brackets.
%
\def\biblabelcontents#1{{\bblrm [#1]}}%
%
% The main text, too, is typeset using \bblrm, which is \rm by default.
%
\def\bblrm{\rm}%
%
% Emphasis for producing, e.g., titles, is done with \it by default.
%
\def\bblem{\it}%
%
% Some styles use a caps-and-small-caps font for author names.  LaTeX
% defines an \sc command but plain TeX doesn't, so we need one here.
% The definition below doesn't load the font unless it's needed, but it
% tries to load only the 10pt version, because it might not exist at
% other point sizes.
%
\def\bblsc{\ifx\@scfont\@undefined
              \font\@scfont = cmcsc10
           \fi
           \@scfont
}%
%
% The major parts of an entry are separated with \bblnewblock.  The
% numbers below are taken from LaTeX's `article' style.
%
\def\bblnewblock{\hskip .11em plus .33em minus .07em }%
%
% Here's where you stick any other bibliography-formatting goodies, or
% redefine the values above.
%
\let\bblhook = \empty
%
%
% Here are the four default definitions for formatting the in-text
% citations.  These are what you redefine (after your \input btxmac but
% before your \bibliography) to get parens instead of brackets, or
% superscripts, or footnotes, or whatever.
%
\def\printcitestart{[}%         left bracket
\def\printcitefinish{]}%        right bracket
\def\printbetweencitations{, }% comma, space
\def\printcitenote#1{, #1}%     comma, space, note (if it exists)
%
% That scheme is pretty flexible.  For example you could use
%     \def\printcitestart{\unskip $^\bgroup}
%     \def\printcitefinish{\egroup$}
%     \def\printbetweencitations{,}
%     \def\printcitenote#1{\hbox{\sevenrm\space (#1)}}
%     \font\eighttt = cmtt8
%     \scriptfont\ttfam = \eighttt
% to get superscripted in-text citations.  (The scriptfont stuff
% exists only to print an undefined citation; it's in cmtt8 because
% there is no cmtt7.)  To get something radically different, however,
% you'll have to define your own \cite command.
%
% When we read `\citation' from the .aux file, it means nothing.
%
\let\citation = \@gobble
%
%
% Now comes the stuff for dealing with LaTeX's \newcommand.  As
% mentioned earlier, this \newcommand will redefine a preexisting
% command; that's different from how LaTeX's \newcommand behaves.
%
\@innernewcount\@numparams
%
% \newcommand{\foo}[n]{text} defines the control sequence \foo to have
% n parameters, and replacement text `text'.
%
\def\newcommand#1{%
   \def\@commandname{#1}%
   \@getoptionalarg\@continuenewcommand
}%
%
% Figure out if this definition has parameters.
%
\def\@continuenewcommand{%
   % If no optional argument, we have zero parameters.  Otherwise, we
   % have that many.
   \@numparams = \ifx\@optionalarg\empty 0\else\@optionalarg \fi \relax
   \@newcommand
}%
%
% \@numparams is how many arguments this command has.  The name of the
% command is \@commandname.  The replacement text for the new macro is #1.
%
\def\@newcommand#1{%
   \def\@startdef{\expandafter\edef\@commandname}%
   \ifnum\@numparams=0
      \let\@paramdef = \empty
   \else
      \ifnum\@numparams>9
         \errmessage{\the\@numparams\space is too many parameters}%
      \else
         \ifnum\@numparams<0
            \errmessage{\the\@numparams\space is too few parameters}%
         \else
            \edef\@paramdef{%
               % This is disgusting, but \loop doesn't work inside \edef,
               % because \body isn't defined.
               \ifcase\@numparams
                  \empty  No arguments.
               \or ####1%
               \or ####1####2%
               \or ####1####2####3%
               \or ####1####2####3####4%
               \or ####1####2####3####4####5%
               \or ####1####2####3####4####5####6%
               \or ####1####2####3####4####5####6####7%
               \or ####1####2####3####4####5####6####7####8%
               \or ####1####2####3####4####5####6####7####8####9%
               \fi
            }%
         \fi
      \fi
   \fi
   \expandafter\@startdef\@paramdef{#1}%
}%
%
%% [[[end of BibTeX-specific stuff]]]
%
%
% Names of references (arguments given in the \cite and \nocite
% commands) and file names (arguments given in the \bibliography and
% \bibliographystyle commands) are recorded in \jobname.aux, called the
% \@auxfile in these macros.  Here's how they get read in.
%
\def\@readauxfile{%
   \if@auxfiledone \else % remember: \@auxfiledonetrue if \noauxfile is defined
      \global\@auxfiledonetrue
      \@testfileexistence{aux}%
      \if@fileexists
         \begingroup
            % Because we might be in horizontal mode when \@readauxfile
            % is called (if it's in response to a \cite or \nocite), we
            % want to ignore all the would-be spaces at the ends of
            % lines in the aux file.  Fortunately, it's highly unlikely
            % an end-of-line might actually be desired.
            % And because we don't change the category code of anything
            % but @, primitives like \gdef can't be used to define labels
            % in the aux file.  The solution adopted by btxmac.tex is to
            % write `\@citedef{LABEL}{DEFINITION}' to the aux file, and
            % use \csname on LABEL.
            \endlinechar = -1
            \catcode`@ = 11
            \input \jobname.aux
         \endgroup
      \else
         \message{\@undefinedmessage}%
         \global\@citewarningfalse
      \fi
      \immediate\openout\@auxfile = \jobname.aux
   \fi
}%
%
% The \@readauxfile macro does all that work the first time it's called.
% Since it's called once for every \cite, \nocite, \bibliography, and
% \bibliographystyle command that the user issues, we need to remember
% whether the work's been done.  It's considered done if we're not to do
% it---that is, if \noauxfile is defined.
%
\newif\if@auxfiledone
\ifx\noauxfile\@undefined \else \@auxfiledonetrue\fi
%
% It's conceivable you'd want to change how other characters are read;
% to do that, change their category code before doing \input btxmac.
%
%
% After reading the .aux file, \@readauxfile opens it for writing.
% The \@writeaux macro does the actual writing (as long as
% \noauxfile is undefined).
%
\@innernewwrite\@auxfile
\def\@writeaux#1{\ifx\noauxfile\@undefined \write\@auxfile{#1}\fi}%
%
%
% A macro package that uses btxmac.tex might define
% \@undefinedmessage (before doing an \input btxmac).
%
\ifx\@undefinedmessage\@undefined
   \def\@undefinedmessage{No .aux file; I won't give you warnings about
                          undefined citations.}%
\fi
%
% Even if citations are undefined, we want to complain only if
% \@citewarningtrue.  The default is to set \@citewarningtrue unless
% \noauxfile is defined.  Again, a macro package that uses
% btxmac.tex might want to redefine this.
%
\@innernewif\if@citewarning
\ifx\noauxfile\@undefined \@citewarningtrue\fi
%
%
% Finally, before leaving we restore @'s old category code.
%
\catcode`@ = \@oldatcatcode

%
% file bib4plain.tex ends  here 
%

\nocite{Al}
%\nocite{Sh:49}
\nocite{Sh:329}
\nocite{GrRoSp}

\topmatter
\title {{\it A finite partition Theorem with double exponential bound} \\
Dedicated to Paul Erd\H{o}s} \endtitle
\rightheadtext {Double Exponential Bound}
\bigskip
\author {Saharon Shelah \thanks{\null\newline 
I thank Alice Leonhardt for excellent typing.\null\newline
Research supported by the United States-Israel Binational Science
Foundation\null\newline 
Publication number 515; Latest Revision 02/95 \null\newline
%Previous Revision 7/10/94
%Previous Revision 5/20/94
Done Erd\H{o}s 80th birthday meeting in Keszthely, Hungary; July 1993}
\endthanks} 
\endauthor
%written 20/7/93
%Previous revision --- 3/94
\affil {Institute of Mathematics \\
The Hebrew University \\
Jerusalem, Israel
\medskip
Rutgers University \\
Department of Mathematics \\
New Brunswick, NJ USA} \endaffil
\abstract
We prove that double exponentiation is an upper bound to Ramsey theorem
for colouring of pairs when we want to predetermine the order of the
differences of successive members of the homogeneous set.
\endabstract
\endtopmatter
\document
\newpage

The following problem was raised by Jouko Vaananen for model theoretic
reasons (having a natural example of the difference between two kinds
of quantifiers, actually his question was a specific case), and propagated by
Joel Spencer: Is there
for any $n$, $c$ an $m$ such that 
\roster
\item "{$(*)$}"  for every colouring $f$ of the pairs from
$\{ 0,1,\dotsc,m-1 \}$ by $2$ (or even $c$) colours, there is a monocromatic
subset $\{ a_0,\dotsc,
a_{n-1} \}$, $a_0 < a_1 < \dots$ such that the sequence
$\langle a_{i+1} - a_i:i < n-1 \rangle$ is with no repetition and is with any
pregiven order.
\endroster
\tracingmacros2
Noga Alon \cite{Al} and independently Janos Pach proved that for every
$n$, $c$ there is such an $m$ as in  $(*)$; 
Alon used van der 
Waerden numbers (see \cite{Sh:329}) (so obtained weak bounds). \newline
Later Alon improved it to iterated exponential (Alan Stacey and also the
author have later
and independently obtained a similar improvement).
  We get a double exponential bound.
The proof continues \cite{Sh:37}.  Within the "realm" of double
exponential in $c,n$ we do not try to save. \newline
We thank Joel Spencer for telling us the problem, and Martin Goldstern
for very careful proof reading.

\demo{Notation}  Let $\ell,k,m,n,c,d$ belong to the set $\Bbb N$ of
natural numbers (which include zero). A sequence $\eta$ is $\langle
\eta(0),\ldots,\eta({\ell g}\eta -1)\rangle$, also $\rho$, $\nu$ are
sequences.
$\eta \triangleleft \nu$ means that $\eta$ is a proper initial
subsequence of $\nu$. 
We consider sequences as graphs of functions (with domain of the form
$n=\{0, \ldots, n-1\}$, so e.g.\ 
$\eta\cap \nu$ means the largest initial segment common to $\eta $ and
$\nu$.     $\eta\extend \langle s\rangle $ is the sequence
$\langle\eta(0),\ldots,\eta({\ell g}\eta -1), s\rangle$ (of length $\ell
g(\eta)+1$). 
 \newline
Let ${}^\ell m := \biggl\{ \eta:{\ell g}(\eta) = \ell,\text{ and Rang}(\eta)
\subseteq \{ 0,\dotsc,m-1 \} \biggr\}$, \newline
${}^{\ell >}m = \dsize \bigcup_{k < \ell} {}^k m$. \newline
For $\nu\in {}^{\ell >}m$ we write $[\nu]_{{}^\ell m}$ (or $[\nu]$ if
$\ell$ and $m$ are clear from the context) for the set 
$\{\eta\in {}^\ell m: \nu \triangleleft \eta\}$. \newline
$[A]^n = \{ w \subseteq A:|w| = n \}$. \newline
Intervals $[a,b),(a,b),[a,b)$ are the usual intervals of integers.
The proof is similar to \cite{Sh:37} but sets are replaced by trees.
\enddemo
\bigskip

\definition{0 Definition}  $r=r(n,c)$ is the first number
$m$ such that:
\medskip
\roster
\item "{$(*)^{n,c}_m$}"  \underbar{for every} $f:[\{ 0,\dotsc,m-1\}]^2
\rightarrow \{ 0,\dotsc,c-1 \}$ and linear order \newline
 $<^*$ on $\{ 0,\dotsc,n-2\}$
\underbar{we can find} $a_0 < \cdots < a_{n} \in [0,m)$ such that:
{\roster
\itemitem{ (a) }  $f \restriction \{ \{ a_i,a_j \}:0 \le i < j < n\}$ is
constant
\itemitem{ (b) }  the numbers $b_\ell := a_{\ell + 1} - a_\ell$ (for
$\ell < n-1)$ are with no repetitions and are ordered by $<^*$, i.e.
$i <^* j \Rightarrow b_i < b_j$ 
\endroster}
\endroster
\enddefinition
\medskip

We will find a double exponential bound for 
$r=r(n,c)$, specifically, 
$r=r(n,c)\le 2^{(c(n+1)^3)^{nc}}$
(so our bound is double exponential in $n$ and in $c$).

This is done in conclusion 5. Alon conjectures that the true order of
magnitude of $r(n,c)$ is single exponential, and Alon and Spencer have
proved this 
for the case where the sequence
 $\langle a_{i+1}-a_i: i< n-1 \rangle$ is monotone.

\bigskip

\definition {1 Definition}  We say $S$ is an $(\ell,m^*,m,u)$-tree if:
\medskip
\roster
\item "{(a)}"  $u \subseteq \{ 0,1,\dotsc,\ell - 1 \}$ and $m^* \ge m$
\item "{(b)}"  $S \subseteq {}^{\ell \ge}(m^*)$
\item "{(c)}"  $S$ is closed under initial segments
\item "{(d)}"  if $\nu \in S$ is $\triangleleft$-maximal \underbar{then} 
$\ell g(\nu) = \ell$
\item "{(e)}"  if $\nu \in S$ and $\ell g(\nu) \in u$ \underbar{then} for
$m$ numbers $j < m^*$ we have $\nu\extend \langle j\rangle \in S$
\endroster
\enddefinition
\bigskip

\proclaim{2 Claim}  Suppose $k,\ell,m,p,m^* \in {\Bbb N}$ satisfies
\medskip

\noindent
$(*)_{k,\ell,m,p,m^*} \qquad \qquad  m^* \ge pm^{\ell k+1}$,
\medskip

\noindent
\underbar{then} for every $i(*) < \ell$, and $u \subseteq [0,\ell-1)$,
and $|u| \le p$
and $f:[{}^\ell (m^*)]^2 \rightarrow \{ 0,1 \}$ \underbar{there is}
$T \subseteq {}^{\ell \ge}(m^*)$ closed under initial segments,
$(T,\triangleleft) \cong ({}^{\ell \ge}m,\triangleleft)$ satisfying
\medskip
\roster
\item "{$\bigoplus$}"  for every $\eta_1,\dotsc,\eta_k \in T \cap
{}^\ell(m^*)$ with $\langle \eta_1 \restriction i(*),\dotsc,\eta_k
\restriction i(*) \rangle$ pairwise distinct we can find a set $S$
such that:
{\roster
\itemitem{ (a) }  $S$ is an $(\ell,m^*,m,u\setminus i(*))$-tree
\itemitem{ (b) }  if $j \in [1,k)$ and $\nu \in S \cap {}^\ell(m^*)$
then $f(\{ \eta_j,\eta_k\}) = f(\{ \eta_j,\nu \})$
\itemitem{ (c) }  $\nu \in S \cap {}^\ell (m^*) \Rightarrow \eta_k
\restriction i(*) \triangleleft \nu$
\endroster}
\endroster
\endproclaim
\bigskip

\remark{Remark}
\item{(1)} With minor change we can demand in $\bigoplus$ ``for
any $i(*)<\ell$''.
\item{(2)} We could use here $f$ with range $\{0, \ldots, c-1\}$,
and in claim 3 get a longer sequence $\langle \nu_\ell: \ell<
n^*\rangle$ such that $f(\{\nu_{\ell_1}, \nu_{\ell_2}\})$ depends just
on $\nu_{\ell_1}\cap \nu_{\ell_2}$, then use a partition theorem on
such colouring.
\endremark

\demo{Proof}  For each $\eta \in {}^{\ell >}(m^*)$ choose randomly a set
$A_\eta \subseteq [0,m^*)$, $|A_\eta|=m$, $A_\eta=\{x^\eta_1,\ldots,
x^\eta_m\}$ (pairwise distinct, chosen by order) (not all are relevant, some
can be fixed). 
\medskip

\noindent We define $T = \{ \eta:\eta \in {}^{\ell \ge}(m^*)$, and $i < \ell
g(\eta) \Rightarrow \eta(i) \in A_{\eta \restriction i} \}$. \newline 
We have a natural isomorphism $h$ from ${}^{\ell \ge}m$ onto $T$:
$$
h(\nu) = \eta \Leftrightarrow \dsize \bigwedge_{i < \ell g \nu}
\eta(i) = x^{\eta \restriction i}_{\nu(i)}.
$$
\noindent
Our problem is to verify $\bigoplus$, we prove that
the probability that it fails
is $< 1$, this suffices. \newline
We can represent it as:
\medskip
\roster
\item "{$(*)_1$}"  if $\nu_1,\dotsc,\nu_k \in {}^\ell m$ and $\nu_1 \restriction
i(*),\dotsc,\nu_k \restriction i(*)$ distinct, then for \newline
$h(\nu_1),\dotsc,h(\nu_k)$ there is $S$ as required there.
\endroster
\medskip

So it suffices to prove that for any given such $i(*) < \ell$ and
$\nu_1,\dotsc,\nu_k$ the probability of failure is
$\displaystyle
< \frac 1{\binom {|{}^\ell m|} k} = \frac 1{\binom {m^\ell} k}$ as it
suffice the demand in $(*)_1$ to hold for the minimal suitable $i(*)$.
Wlog $u=u\setminus i(*)$. \newline
For this we can assume $x^\rho_j$ are fixed whenever $\neg[h(\nu_k)
\restriction i(*)) \trianglelefteq \rho]$ or $\ell g(\rho) \notin u$. \newline
Let $Y = \{ \eta \in {}^\ell(m^*)$:Prob$[h(\nu_k) = \eta] \ne 0\}$, so
$|Y| = m^{|u|}$. \newline
So $h(\nu_1),\dotsc,h(\nu_{k-1})$ are determined.  Now
$h(\nu_1),\dotsc,h(\nu_{k-1})$ and $f$ induces an equivalence relation $E$ on
$Y$:
$$
\eta'E \eta'' \text{ iff } \dsize \bigwedge^{k-1}_{j=1} f(\{ h(\nu_j),\eta'\})
= f(\{h(\nu_j),\eta'' \}).
$$
The number of classes is $\le 2^{k-1}$, let them be $A_1,\dotsc,
A_{2^{k-1}}$ (they are pairwise disjoint, some may be empty).

We call $A_j$ \underbar{large} if there is $S$ as required in clauses
(a) and (c) of
$\bigoplus$ such that \newline
$(\forall \rho)[\rho \in S \cap {}^\ell(m^*)
\Rightarrow \rho \in A_j]$.
\medskip

\noindent
It is enough to show that the probability of $h(\nu_k)$ belonging to a
non-large 
equivalence class is $< \frac 1{\binom {m^\ell} k}$, hence it is 
enough to prove:
\medskip
\roster
\item "{$(*)_2$}"   $A_j$ not large $\Rightarrow$ Prob$(h(\nu_k) \in A_j) <
\frac 1{\binom {m^\ell} {k} \times 2^k}$.
\endroster
\medskip

\noindent
So assume $A_j$ is not large.  Let $Y^* := \{ \eta \restriction i:\eta \in Y
\text{ and }i \le \ell \}$.
$$\align
\text{Let } Z_j := \biggl\{ \eta \in Y^*\;:\;&\text{there is } S \subseteq Y^*
\text{ satisfying} \\
  &(a)' \quad S \text{ is an } (\ell,m^*,m,u \backslash \ell g(\eta))
\text{-tree} \\
  &(b)' \quad \nu \in S \cap {}^\ell(m^*) \Rightarrow \eta \trianglelefteq
\nu \\
  &(c)' \quad \text{for every } \nu \in S \cap {}^\ell(m^*) \text{ we have }
\nu \in A_j 
	\smash{\biggr\}.}
\endalign
$$
Let $Z^*_j := \{ \eta \in Z_j$: there is no
$\eta' \triangleleft \eta,\eta' \in  Z_j \}$. \newline
Clearly $h(\nu_k \restriction i(*)) \notin Z_j$, (as $A_j$ is not large)
hence $h(\nu_k \restriction i(*)) \notin Z^*_j$. \newline
Clearly
\medskip
\roster
\item "{$(*)_3$}"  for $\eta \in Y^*\setminus Z_j$ such that
${\ell g}(\eta)\in u$ we have
$$|\{ i:\eta \extend\langle i\rangle \in Z^*_j\text{ (or even }\in
Z_j) \}| < m.$$ 
\endroster
But if $\nu\triangleleft\eta$, $\nu\in Z_j$ then $\eta\notin Z^*_j$. Hence
\roster
\item "{$(*)_4$}"  for $\eta \in Y^*$ such that ${\ell g}(\eta)\in u$ we have
$$|\{ i:\eta \extend\langle i\rangle \in Z^*_j\}| < m.$$ 
\endroster

\noindent Now 
\roster
\item "{$(*)_5$}"  if $\eta \in A_j$ (hence $\eta \in Y$; remember that $A_j$
is not large) \underbar{then} $\eta \in Z_j,$ \newline
hence $\dsize \bigvee_{j\in u} \eta \restriction (j+1) \in Z^*_j$.
\endroster
\medskip

\noindent
So
$$
\align
\text{Prob}(h(\nu_k) \in A_j) &\le \text{ Prob}
\biggl( \dsize \bigvee_{j \in u} [h(\nu_k \restriction (j+1)) \in Z^*_j]
\biggr) \\
  &\le \dsize \sum_{j \in u} \text{ Prob}(h(\nu_k \restriction (j+1) \in
Z^*_j)) \\
  &< |u| \times {\frac m{m^*}}.
\endalign
$$
\noindent
(first inequality by $(*)_5$, second inequality trivial, last inequality
by $(*)_4$ above).\newline
So it suffices to show:
$$
|u| \times {\frac m{m^*}} \le \frac 1{\binom {m^\ell}{k} \times 2^k}
$$
equivalently
$$
m^* \ge |u| \times m \times {\binom {m^\ell} k} \times 2^k
$$
as ${\binom {m^\ell} k} \le m^{\ell k}/k!$, and $|u| \le p$ by the hypothesis
$(*)_{k,\ell,m,p,m^*}$ we finish.\hfill$\square_2$
\enddemo
\bigskip

\proclaim{3 Lemma}  Assume
\medskip
\roster
\item "{(a)}"  $\rho_1,\dotsc,\rho_n \in {}^{n-1}2$ are distinct,
for $\ell\in \{2, \ldots, n\}$, we have $r_\ell \in \{1, \ldots, \ell-1\}$
such that $\ell g(\rho_\ell \cap \rho_{r_\ell})=\ell-1$ and $r\in \{1, \ldots,
\ell-1\} \setminus \{r_\ell\} \Rightarrow \ell g(\rho_\ell \cap
\rho_{r}) < \ell$
\item "{(b)}"  $f:[{}^\ell m]^2 \rightarrow [0,c) = \{ 0,\dotsc,c-1 \}$
\item "{(c)}"  $m = 2^{(n+1)^{(c+1)n}}$.
\endroster
\medskip

\noindent
\underbar{Then} we can find $\eta_1,\dotsc,\eta_n \in {}^\ell m$ such that:
\medskip
\roster
\item "{$(\alpha)$}"  $f \restriction [\{ \eta_1,\dotsc,\eta_n \}]^2$
is a constant function
\item "{$(\beta)$}"  $\langle \ell g(\eta_{i+1} \cap \eta_{r_{i+1}}):
i= 1,\dotsc, n-1\} \rangle$ is a sequence with no repetitions ordered just like
$\langle \ell g(\rho_{i+1} \cap \rho_{r_{i+1}}): i = 1,n-1\} \rangle$; also:
\endroster
$$
\align
\eta_{i+1}(\ell g(\eta_{i+1} \cap \eta_{r_{i+1}})) &
<\eta_{r_{i+1}}(\ell g(\eta_i \cap
\eta_{r_{i+1}})) \\
  &\Leftrightarrow \rho_{i+1}(\ell g(\rho_{i+1} \cap \rho_{r_{i+1}})) <
\rho_{r_{i+1}}(\ell g(\rho_{i+1} \cap \rho_{r_{i+1}})).
\endalign
$$
\endproclaim

\remark{3A Remark} {\it 1)}\ \ Note that if $\Gamma\subseteq {}^{n-1}2$,
$|\Gamma|=n$ and the set $\{\rho_1\cap\rho_2: \rho_1\neq\rho_2$ are from
$\Gamma\}$ has no two distinct members with the same length then we can list
$\Gamma$ as $\langle\rho_1,\ldots,\rho_n\rangle$ as required in clause (a) of
lemma 3.

\noindent{\it 2)}\ \ So if $<^*$ is a linear order on $\{1,\ldots,n-1\}$ then
we can find distinct $\rho_1,\ldots,\rho_n\in {}^{n-1}2$ as in clause (a) of
lemma 3 and a permutation $\sigma$ of $\{1,\ldots,n\}$ such that:

for $i\neq j\in \{1,\ldots,n-2\}$ we have
$$
i<^*j\quad \text{iff}\quad {\ell g}(\rho_{\sigma(i)}\cap \rho_{\sigma(i+1)})
> {\ell g}(\rho_{\sigma(j)}\cap \rho_{\sigma(j+1)}).
$$
(E.g. use induction on $n$.)
\endremark

\demo{Proof}  Let us define $\langle m_j:2 \le j \le cn \rangle$ by
induction on $j$:
$$m_2 = n^c,\qquad m_{j+1} = n^c m_j^{n^c(n+1)+1}.
$$
Check that $m_j \le 2^{(n+1)^{(c+1)j}} $, so in particular
 $ m_{(cn)} \le m$.  Now we claim that for any number $d \in [1,c]$ the
following holds:
\medskip
\roster
\item "{$\bigotimes_d$}"  Assume $q \in [0,n^c)$ such that $q$ divisible
by $n^d$ and $q + n^d \le cn$,
$u = [q,q + n^d)$, $T^*$ is an $(\ell,m,m_{(q+dn)},u)$-tree and $f$ is
a function 
from $[T^* \cap {}^\ell m]^2$ with range of cardinality $d$. \underbar{Then}
we can
find $\eta_1,\dotsc,\eta_n \in T^* \cap ({}^\ell m)$ such that clauses
$(\alpha)$ and $(\beta)$ of the conclusion of lemma 3 hold.
\endroster
\medskip

\noindent
This suffices: use $q=0, d=c$.  We prove this by induction on $d$.  If
%% $d=2$, trivial as $m_j = 2$ as only one colour occurs.  For $d+1$,
%% without
$d=1$, trivial as only one colour occurs.  For $d+1 > 1$, without
loss of generality $\text{Rang}(f) = [0,d]$, let $f':[{}^\ell m]^2
\rightarrow \{ 0,1 \}$ be $f'(\{ \eta',\eta''\}) = \text{ Min}\{f (
\{\eta',\eta''\},1\}$. Let for $j < n$, $u_j := [q+n^d j,q+n^d j + n^d)$.
By downward induction on $j \in [1,n]$ we try to define $T_j$ such that:
\medskip
\roster
%\widestnumber\item{(iii)}
\item "{(i)}"  $T_j$ is an $(\ell,m,m_{(q+nd+j)},
\dsize \bigcup_{i < j} u_i)$-tree
\item "{(ii)}"  $T_j \subseteq T_{j+1}$
\item "{(iii)}"  for every $j \in [1,n-1]$ and $\eta_1,\dotsc,\eta_j \in
T_j \cap {}^\ell m$ with $\langle\eta_1\restriction (q+n^dj),
\dotsc,\eta_j \restriction (q+n^dj) \rangle$ pairwise distinct, we
can find $\eta' \ne \eta'' \in T_{j+1} \cap {}^\ell m$ such that:
\endroster

$$
\eta_j \restriction (q+n^dj) \triangleleft \eta',\eta''
$$
$$
\ell g(\eta' \cap \eta'' ) \in u_j
$$
$$
f'(\{ \eta',\eta'' \}) = 0
$$
$$
\dsize \bigwedge_{t \in [1,j)} \biggl[ 0 = f'( \{ \eta_t,\eta_j \})
\Rightarrow
0= f'(\{ \eta_t,\eta'\}) = f'(\{ \eta_t,\eta''\}) \biggr].
$$

\noindent 
\underbar{This suffices} as then we can choose by induction on
$j=1,\ldots,n$ a sequence $\nu_j \in T_j \cap ({}^\ell m)$
such that (after reordering) the set $\{\nu_1, \ldots, \nu_n\}$ will
serve as $\{\eta_1,\ldots, \eta_n\}$ of $\otimes_d$ (with the constant
colour being zero). Let us do it in detail. 

By induction on $j=1, \ldots, n$ we choose $\nu^j_1, \ldots,
\nu^j_{j-1}$ such that:
\item {(a)} $\nu^j_1, \ldots, \nu^j_j$ are distinct members of
$T_j\cap ({}^\ell m)$
\item {(b)} $f\restriction [\{\nu^j_1, \ldots, \nu^j_j\}]^2$ is
constantly zero
\item {(c)} for $\ell = \{2, \ldots, j\}$ we have $\ell
g({\nu^j_\ell} \cap {\nu^j_{r_\ell}})\in u_{\ell-1}$

\noindent For $j=1$ no problem. In the induction step, i.e. for $j+1$, we
apply the condition (iii) above with $\langle\nu^j_\ell: \ell\in [1,j]$, 
$\ell \neq q_{j+1}\rangle \extend
\langle \nu^j_{q_{j+1}}\rangle$ here standing for $\eta_1, \ldots,
\eta_j$ there (we want $\nu^j_{q_{j+1}}$ be the last), the condition 
$\langle \eta_\ell \restriction (q+ n^dj) : \ell=1, \ldots, j\rangle$
with no repetition follows by clause (c), so we get $\eta'$, $\eta''
\in T_{j+1} \cap ({}^\ell m)$ as there. W.l.o.g. $\eta'(\ell
g(\eta'\cap \eta''))< \eta''(\ell g(\eta'\cap \eta''))$.

We now define $\nu^{j+1}_\ell$ for $\ell=1, \ldots, j+1$:
\itemitem{ } If $\ell\in \{1, \ldots, j+1\} \setminus \{j+1,
r_{j+1}\}$ then $\nu^{j+1}_\ell= \nu^j_\ell$ (remember $T_j \subseteq
T_{j+1}$).
\itemitem{ } If $\rho_{j+1}(\ell g(\rho_{j+1}\cap \rho_{r_{j+1}})) <
\rho_{r_{j+1}}(\ell g(\rho_{j+1} \cap \rho_{r_{j+1}}))$
\itemitem{ }\ \ \ \ then $\nu^{j+1}_{j+1}=\eta'$ and
$\nu^{j+1}_{r_{j+1}}=\eta''$. 
\itemitem{ } If $\rho_{r_{j+1}}(\ell g(\rho_{j+1} \cap \rho_{r_{j+1}})) <
\rho_{j+1} (\ell g(\rho_{j+1} \cap \rho_{r_{j+1}}))$
\itemitem{ }\ \ \ \ then $\nu^{j+1}_{j+1}=\eta''$ and
$\nu^{j+1}_{r_{j+1}}=\eta'$. 
\newline

Now check. (Note $T_0$, $u_0$ could be omitted above.)
\enddemo
\bigskip

\definition{Carrying the Inductive Definition}  For $j = n$, trivial: let
$T_n = T^*$ (given in $\otimes_d$). \newline
For $j \in [2,n)$, where $T_{j+1},\dotsc,T_n$ are already defined, we apply
Claim 2 with
$$
m_{q+nd+j+1},\ m_{q+nd+j},\ n^d(j{+}1),\ j,\ q{+}n^dj,\
[q{+}n^dj,q{+}n^dj{+}n^d),\ n^dj
$$
here\ \ \ \ standing for
$$
m^*,\ m,\ \ell,\ k,\ i(*),\ u,\ p
$$
there. (I.e. the tree in claim 2 is replaced by one isomorphic to it,
levels outside $\bigcup_{i<j} u_i$ can be ignored.)

\noindent So we need to check $m_{q+nd+j+1} \ge n^dj
(m_{q+nd+j})^{n^d(j+1)j+1}$, \
which holds by the definition of the $m_i$'s (as $d\leq c_1$).
\underbar{But} we require above more than in Claim 2
(preferring the colour $0$).  But if it fails for $T_j$ then for some
$\eta_1,\dotsc,\eta_j$ in $T_j$ we have $S$ as in $\oplus$ of Claim 2,
with no $\eta^1 \ne \eta^2 \in S \cap {}^\ell m$ such that $f'(\{ \eta',\eta''
 \}) = 0$.  On $S$ we can apply our induction hypothesis on $d$ --- allowed as
the original $f$ misses a colour (the colour zero) when restricted to
$S$.\hfill$\square_3$ 
\enddefinition
\bigskip

\proclaim{4 Fact} 
Let $\langle \eta_i: i < (2m-1)^\ell \rangle $ enumerate ${}^\ell (2m-1)$
in lexicographic order. 
Let
$$
A=\{0, 2, \ldots, 2m-2\}\subseteq [0, 2m-1),\text{ so }|A|=m,
$$
$$
B:= \{ i < (2m-1)^\ell: \eta_i \in {}^\ell A \}.
$$
Let
$$
{\IT low}_{m,\ell}(k) = (2m-1)^{\ell-k-1},\text{ and }
$$
$$
{\IT high}_{m,\ell}(k)=  (2m-1)^{\ell-k}-1.
$$
Then:
\roster
\item"{(0)}" if $i\neq j$ are in $B$, $|\eta_i\cap \eta_j|=k$ then:
$$
i<j \text{ iff } \eta_i(k)< \eta_j(k).
$$
\item"{(1)}"  If $i<j$ are in $B$, $|\eta_i\cap \eta_j|=k$, then:
$$ 
{\IT low}_{m,\ell}(k)\  \le\  j-i\  \le\ {\IT high}_{m,\ell}(k)
$$
\item"{(2)}"  If $i_1 < j_1$, $i_2 < j_2$ are all in $B$, and 
	$|\eta_{i_1}\cap \eta_{j_1} | \not= 
	|\eta_{i_2}\cap \eta_{j_2} |$,
then: 
$$ 
 j_1-i_1 < j_2 -i_2 \text { \ \  iff \ \  } 
|\eta_{i_1}\cap \eta_{j_1} | > |\eta_{i_2}\cap \eta_{j_2} |.$$
\endroster
\endproclaim

\demo{Proof} (0) Check.
\newline
(1) Let $\nu:= \eta_i\cap \eta_j$ and $k:= |\nu|$.   We are
looking for upper and lower bounds of the cardinality of the set (the
order is lexicographic)
$$ C:= \{ \eta\in {}^\ell(2m-1): \eta_i \le \eta < \eta_j\}.$$
Clearly each $\eta\in C $ satisfies $\nu \triangleleft \eta$.  Moreover, 
since $\eta_j\in {}^\ell m $, each element $\eta\in C$ must satisfy
$\eta(k) \le \eta_j(k) < 2m-1$. Hence
$$ C \subseteq \bigcup_{s<2m-1 } [\nu\extend \langle
s\rangle]_{({}^\ell(2m-1))}\setminus \{\eta_j\} $$
so we get  $|C| \le  (2m-1)\cdot (2m-1)^{\ell-k-1}-1=(2m-1)^{\ell-k}-1$. \newline
For $k=\ell-1$ the lower bound claimed in (1) is trivial, so assume
$k\le \ell- 2$.  Let $\nu':= (\eta_i\restriction k)\extend \langle
\eta_i(k)+1\rangle$ (note: as $\eta_i(k) < \eta_j(k)$ are both in $A$,
necessarily $\nu'(k)<\eta_j(k)$).  Then we have 
$$ C \supseteq \bigcup_{s < 2m-1} [\nu'\extend \langle
s\rangle]_{({}^\ell(2m-1))}\cup \{ \eta_i\},$$
so $|C| \ge (2m-1)^{\ell-k-1}+1$.

\noindent
Proof of (2):  Check that  ${\IT high}_{m, \ell}(k)<{\IT low}_{m, \ell}(k+1)$. 
\enddemo\hfill$\square_4$

\remark{Remark} Also ${\IT low}_{m,\ell}(k)=(2m-1)^{\ell-k-1}+1$ is O.K. but
with the present bound we can use only $\eta_i$ with $\eta_i(\ell-1)<m$,
$B=\{i: \eta_i\restriction (\ell-1)\in {}^\ell A\}$. So $(2m-1)^\ell$ can be
replaced by $m(2m-1)^{\ell-1}$.
\endremark

\proclaim{5 Conclusion}  If $f:[0,(2m-1)^\ell)^{[2]} \rightarrow
[0,c)$, $\ell=nc$, $m := 2^{(n+1)^{(c+1)n}}$, then we can find
$a_0 < \dots < a_{n-1} < m$ such that $f \restriction [\{ a_1,\dotsc,
a_{n-1} \}]^2$ constant and
$\langle a_{i+1} - a_i:i < n - 1 \rangle$
in any pregiven order
\endproclaim

\demo{Proof} As in fact 4,
let $\langle \eta_i: i < (2m-1)^\ell
\rangle $ enumerate ${}^\ell (2m-1)$ in lexicographic order, and let
 $B:= \{ i < (2m-1)^\ell: \eta_i \in {}^\ell m \}$ and for $i\in B$
let $\eta'_i= \langle \eta_i(\ell)2: i<\ell\rangle$.   Define a
function 
$f':[{}^\ell m]^2 \to c$ by requiring $f'(\{\eta'_i, \eta'_j\}) =
f(\{i,j\})$ for all $i,j\in B$. 
Now the conclusion follows from lemma 3, Remark 3A(2) and Fact 4, particularly
clause (2). 

% 
% 
%  Translate to Lemma 3.  I.e. the set ${}^\ell m$ has $m^\ell$
% elements and it is linearly ordered by $<_{\ell x}$, the lexicographic order;
% so let $\{ \eta_i:i < m^\ell \}$ enumerate ${}^\ell m$ such that
% $i < j < m^\ell \Rightarrow \eta_i <_{\ell x} \eta_j$.  Define a function
% $f':[{}^\ell m]^2 \rightarrow [0,c)$ by $f'(\{ \eta_i,\eta_j \}) =
% f(\{ i,j \})$.  Now the conclusion of Lemma 3 gives the desired conclusion
% because:
% \medskip
% \roster
% \item "{$(*)$}"  if $i_1 < j_1 < m^\ell,i_2 < j_2 < m^\ell$ and
% $\ell g(\eta_{i_1},n_{j_1}) \ne \ell g(\eta_{i_2} \cap \eta_{j_2})$
% \underbar{then} $j_1 - i_1 < j_2 - i_2 \Leftrightarrow \ell g(\eta_{i_1} \cap
% \eta_{j_1}) > \ell g(\eta_{i_2} \cap \eta_{j_2})$.\hfill$\square_4$
% \endroster

\enddemo\hfill$\square_5$
\newpage

\head{REFERENCES}\endhead
\bigskip

\bibliography{listb,listx}
\bibliographystyle{literal-unsrt}

\shlhetal
\enddocument
\bye